\documentclass[12pt]{article}
\usepackage{amsfonts,amssymb,epsfig}
\usepackage{amsmath}

\date{}
\begin{document}
\newtheorem{df}{Definition}
\newtheorem{thm}{Theorem}

\newtheorem{lm}{Lemma}
\newtheorem{pr}{Proposition}
\newtheorem{co}{Corollary}
\newtheorem{re}{Remark}
\newtheorem{note}{Note}
\newtheorem{claim}{Claim}
\newtheorem{problem}{Problem}

\def\R{{\mathbb R}}

\def\E{\mathbb{E}}
\def\calF{{\cal F}}
\def\N{\mathbb{N}}
\def\calN{{\cal N}}
\def\calH{{\cal H}}
\def\n{\nu}
\def\a{\alpha}
\def\d{\delta}
\def\t{\theta}
\def\e{\varepsilon}
\def\t{\theta}
\def\g{\gamma}
\def\G{\Gamma}
\def\b{\beta}
\def\pf{ \noindent {\bf Proof: \  }}

\newcommand{\qed}{\hfill\vrule height6pt
width6pt depth0pt}
\def\endpf{\qed \medskip} \def\colon{{:}\;}
\setcounter{footnote}{0}

\renewcommand{\qed}{\hfill\vrule height6pt  width6pt depth0pt}

\title{Injective Tauberian operators on $L_1$ and operators with dense range on $\ell_\infty$
}

\author{William B. Johnson\thanks{Supported in part by NSF DMS-1301604  and U.S.-Israel Binational Science
Foundation
 }, Amir Bahman Nasseri,\\   Gideon Schechtman\thanks{Supported in part by U.S.-Israel Binational Science Foundation. Participant
 NSF Workshop in Analysis and Probability, Texas A\&M University
 }, \ and Tomasz Tkocz\thanks{T. Tkocz thanks his PhD supervisor, Keith Ball, for his invaluable constant
advice and encouragement}} \maketitle

\begin{abstract}
There exist injective Tauberian operators on $L_1(0,1)$ that have dense, non closed range. This gives injective, non surjective operators on $\ell_\infty$ that have dense range. Consequently, there are two quasi-complementary, non complementary  subspaces of $\ell_\infty$ that are isometric to $\ell_\infty$.

\end{abstract}

\section{Introduction}

A (bounded, linear) operator $T$ from a Banach space $X$ into a Banach space $Y$ is called Tauberian provided $T^{**-1}Y=X$.  The structure of Tauberian operators when the domain is an $L_1$ space is well understood and exposed in Gonz\'ales and Mart\'inez-Abej\'on's book \cite[Chapter 4]{gon}. (For convenience they only consider $L_1(\mu)$ when $\mu$ is finite and purely nonatomic, but their proofs for the results we mention work for general $L_1$ spaces.) In particular,
\cite[Theorem 4.1.3]{gon}    implies that when $X$ is an $L_1$ space,   an  operator $T:X\to Y$ is Tauberian iff whenever $(x_n)$ is a sequence of disjoint unit vectors, there is an $N$ so that the restriction of $T$ to $[x_N]_{n=N}^\infty$ is an isomorphism (and, moreover, the norm of the inverse of the restricted operator is bounded independently of the disjoint sequence). From this it follows that an injective operator $T:X\to Y$ is Tauberian iff it isomorphically preserves isometric copies of $\ell_1$ in the sense that  the restriction of $T$  to any subspace of $X$ that is isometrically isomorphic to $\ell_1$ is an isomorphism. (Recall that a subspace of an $L_1$ space is isometrically isomorphic to $\ell_1$ iff it is the closed linear span of a sequence of non zero disjoint vectors  \cite[Chapter 14.5]{roy}.) Since $Tu$ is Tauberian if $T$ is Tauberian and $u$ is an isomorphism, one deduces that an injective Tauberian operator from an $L_1$ space isomorphically preserves isomorphic copies of $\ell_1$ in the sense that  the restriction of $T$ to any subspace of $X$ that is isomorphic to $\ell_1$ is an isomorphism. Thus injective Tauberian operators from an $L_1$ space are opposite to $\ell_1$-singular operators; i.e., operators whose restriction to every  subspace isomorphic to $\ell_1$ is {\em not} an isomorphism.

The main result in this paper is a negative solution to \cite[Problem 1]{gon}: Suppose $T$ is a Tauberian operator on an $L_1$ space.  Must $T$ be upper semi-Fredholm; i.e., must the range $\mathcal{R}(T)$  of $T$ be closed and the null space $\mathcal{N}(T)$ of $T$ be finite dimensional?  The basic example is a Tauberian operator on $L_1(0,1)$ that has infinite dimensional null space.  This is rather striking because the Tauberian condition is equivalent to the statement that there is $c>0$ so that the restriction of the operator to $L_1(A)$ is an isomorphism whenever the subset $A$ of $[0,1]$ has Lebesgue measure at most $c$.

In fact, we show that there is an injective, dense range, non surjective Tauberian operator on $L_1(0,1)$.  Since $T$ is Tauberian, $T^{**}$ is also injective, so  $\mathcal{R}(T^*)$ is dense and proper, and  $T^*$ is injective because  $\mathcal{R}(T)$ is dense. This solves a problem \cite{bah} the second author raised on MathOverFlow.net that led to the collaboration of the authors.

\section{The examples}

We begin with a lemma that is an easy consequence of characterizations of Tauberian operators on $L_1$ spaces.

\begin{lm}\label{finite} Let $X$ be an $L_1$ space and  $T$ an operator from $X$ to a Banach space $Y$.  The operator $T$ is Tauberian if and only if there is $r>0$ and a natural number $N$ so that if $(x_n)_{n=1}^N$ are disjoint unit vectors in $X$, then $\max_{1\le n \le N} \|Tx_n\| \ge r$.
\end{lm}
\pf The condition in the lemma clearly implies that if $(x_n)$ is a disjoint sequence of unit vectors in $X$, then $\liminf_n \|Tx_n\| >0$, which is one of the equivalent conditions for $T$ to be Tauberian \cite[Theorem 4.1.3]{gon}.  On the other hand, suppose that there are disjoint collections $(x^n_k)_{k=1}^n$, $n=1,2,\dots$ with $ \max_{1\le k \le n } \|T{x^n_k} \| \to 0$ as $n \to \infty$.  Then the closed sublattice generated by $\cup_{n=1}^\infty (x^n_k)_{k=1}^n$ is a separable abstract $L_1$ space (meaning that it is a Banach lattice such that $\|x+y\|= \|x\|+\|y\|$ whenever $|x| \vee  |y| =0$)  and hence is  order isometric to $L_1(\mu)$ for some probability $\mu$ by Kakutani's theorem (see e.g. \cite[Theorem 1.b.2]{lt2}). Choose $1\le k(n) \le n$ so that the support of $x_{k(n)}^n$ in $L_1(\mu)$ has measure at most $1/n$.  Since $T$ is Tauberian, by \cite[Proposition 4.1.8]{gon} necessarily $\liminf_n \|Tx_{k(n)}^n \| >0$, which is a contradiction. \qed

\medskip
The reason that Lemma \ref{finite} is useful for us is that the condition in the Lemma is stable under ultraproducts.  Call an operator that satisfies the condition in Lemma \ref{finite} $(r,N)$-Tauberian.  For background on ultraproducts   of Banach spaces and of operators, see \cite[Chapter 8]{djt}.  We use the fact that the ultraproduct of $L_1$ spaces is an abstract $L_1$ space and hence   is order isometric to $L_1(\mu)$ for some measure $\mu$.

\begin{lm}\label{ultra} Let $(X_k)$ be a sequence of $L_1$ spaces, and for each $k$  let $T_k$  be a norm one linear operator from $X_k$ into a Banach space $Y_k$.  Assume that   there is $r>0$ and a natural number $N$ so that each operator $T_k$ is  $(r,N)$-Tauberian. Let $\mathcal{U}$ be a free ultrafilter on the natural numbers.  Then $(T_k)_\mathcal{U}: (X_k)_\mathcal{U} \to (Y_k)_\mathcal{U}$
is $(r,N)$-Tauberian.
\end{lm}
Here $(T_k)_\mathcal{U}$ is the usual ultraproduct of the sequence $(T_k)$, defined by  \hfill\break

$(T_k)_\mathcal{U}(x_k) = (T_kx_k)$.

\vskip9pt
\pf  The vectors $(x_k)$ and $(y_k)$ are disjoint in the abstract $L_1$ space $(X_k)_\mathcal{U}$ iff $\lim_\mathcal{U} \||x_k|\wedge |y_k|\| =0$, so it is only a matter of proving that if $T$ is $(r,N)$-Tauberian from some $L_1$ space $X$, then for each $\e>0$ there is $\d>0$ so that if $x_1,\dots, x_N$ are unit vectors in $X$ and $\| |x_n|\wedge |x_m| \| < \d $ for $1\le n<m \le N$, then $\max_{1\le n \le N }\|Tx_n\| > r-\e$. But if $x_1,\dots, x_N$ are unit vectors that are $\e$-disjoint as above, and $y_1,\dots,y_n$ are defined by \hfill\break

$y_n := [|x_n|-(|x_n|\wedge (\vee \{|x_m| : m\not= n\})]\text{sign}(x_n)$, \hfill\break

\noindent then the $y_n$ are disjoint and all have norm at least $1-N\d$.  Normalize the $y_n$ and apply the $(r,N)$-Tauberian condition to this normalized disjoint sequence to see that $\max_{1\le n \le N }\|Tx_n\| > r-\e$ if $\d = \d(\e, N)$ is sufficiently small. \qed

\medskip
An example that answers   \cite[Problem 1]{gon} is the restriction of an ultraproduct of operators on finite dimensional $L_1$ spaces constructed in \cite{bgiks}.

\begin{thm}\label{theorem1}
There is a Tauberian operator $T$ on $L_1(0,1)$ that has an infinite dimensional null space.  Consequently, $T$ is not upper semi-Fredholm.
\end{thm}
\pf An immediate consequence of \cite[Proposition 6 \& Theorem 1]{bgiks} is that there is $r>0$ and a natural number $N$ so that for all sufficiently large $n$ there is a norm one   $(r,N)$-Tauberian operator $T_n$ from $\ell_1^n$ into itself with $\dim \mathcal{N}(T_n) > rn$. The ultraproduct $\tilde{T}:=(T_n)_\mathcal{U}$ is then a norm one
$(r,N)$-Tauberian operator on the gigantic $L_1$ space $X_1:=(\ell_1^n)_\mathcal{U}$, and the null space of $\tilde{T}$ is infinite dimensional.  Take any separable infinite dimensional subspace $X_0$ of $\mathcal{N}(\tilde{T})$ and let $X$ be the closed sublattice of $X_1$   generated by $X_0$.  Let $Y$ be the sublattice of $X_1$ generated by $\tilde{T}X$ and let $T$ be the restriction of $\tilde{T}$ to $X$, considered as an operator into $Y$.  So $X$ and $Y$ are separable $L_1$ spaces and by Lemmas \ref{finite} and \ref{ultra} the operator $T$ is Tauberian.  Of course, by construction $\mathcal{N}(T)$ is infinite dimensional and reflexive (because $T$ is Tauberian).  Thus $X$ is not isomorphic to $\ell_1$ and hence is isomorphic to $L_1(0,1)$.  So is $Y$, but that does not matter: $Y$, being a separable $L_1$ space, embeds isometrically into $L_1(0,1)$. \qed

We want to ``soup up" the operator $T$ in Theorem \ref{theorem1} to get an injective, non surjective, dense range Tauberian operator on $L_1(0,1)$.  We could quote a general result \cite[Theorem 3.4]{go} of Gonz\'alez and Onieva to shorten the presentation, but we prefer to give a short direct proof.

We recall a simple known lemma:

\begin{lm}\label{nuclear}
Let $X$ and $Y$ be separable infinite dimensional Banach spaces and $\e>0$. Let $Y_0$ be a countable dimensional dense subspace of $Y$.Then there is a nuclear operator $u:X\to Y$ so that $u$ is injective and $\|u\|_\wedge < \e$ and $uX \supset Y_0$.
\end{lm}
\pf Recall that an $M$-basis for a Banach space $X$  is a biorthogonal system $(x_\a,x^*_\a)\subset X\times X^*$ such that the linear span of $(x_\a)$ is dense in $X$ and $\cap_\a \mathcal{N}(x^*_\a) =\{0\}$.  Every separable Banach space $X$ has an $M$-basis \cite{mac}; moreover, the vectors $(x_\a)$  in the $M$-basis can span any given countable dimensional dense subspace of $X$.

Take $M$-bases $(x_n,x^*_n)$ and $(y_n,y^*_n)$ for $X$ and $Y$, respectively, normalized so that $\|x_n^*\|=1=\|y_n\|$ and such that the linear span of $(y_n)$ is $Y_0$. Choose $\lambda_n>0$ so that $\sum_n \lambda_n < \e$ and set $u(x)= \sum_n \lambda_n \langle x^*_n, x \rangle y_n$. \qed

\begin{thm}\label{theorem2}
There is an injective, non surjective, dense range Tauberian operator on $L_1(0,1)$.
\end{thm}
\pf By Theorem \ref{theorem1} there is a  Tauberian operator $T$ on $L_1(0,1)$ that has an infinite dimensional null space.  By Lemma \ref{nuclear} there is a nuclear operator $\tilde{v} : \mathcal{N}(T) \to L_1(0,1)$ that is injective and has dense range, and we can extend $\tilde{v} $  to a nuclear operator $v$  on $L_1(0,1)$. We can choose $\tilde{v}$ so that $\tilde{v}(\mathcal{N}(T) ) \cap TL_1(0,1)$ is infinite dimensional by the last  statement in Lemma \ref{nuclear}.  This guarantees that the Tauberian operator $T_1:= T+v$ has an infinite dimensional null space (this allows us to avoid breaking the following argument into cases).

Now $\mathcal{N}(T_1) \cap \mathcal{N}(T) = \{0\}$, so again by Lemma \ref{nuclear} and the extension property of nuclear operators there is a nuclear operator $u : L_1(0,1)/{\mathcal{N}(T)} \to \ell_1$ so that the restriction of $u$ to $ Q_{\mathcal{N}(T)} \mathcal{N}(T_1)$ is injective and has dense range (here for a subspace $E$ of $X$, the operator $Q_E$ is the quotient mapping from $X$ onto $X/E$). Finally, define $T_2 : L_1(0,1) \to L_1(0,1) \oplus_1 \ell_1$ by $T_2 x := T_1x \oplus u Q_{\mathcal{N}(T)} x$.  Then $T_2$ is an injective Tauberian operator with dense range.  $T_2$ is not surjective because $P_{\ell_1} T_2$ is nuclear by construction, where $P_{\ell_1}$ is the   projection of $L_1(0,1) \oplus_1 \ell_1$ onto $\{0\} \oplus_1 \ell_1$.  Since $L_1(0,1) \oplus_1 \ell_1$ is isomorphic to $L_1(0,1)$, this completes the proof. \qed

\begin{co} \label{corollary} There is an injective, dense range, non surjective operator on $\ell_\infty$.  Consequently, there is a quasi-complementary, non complementary decomposition of $\ell_\infty$ into two subspaces each of which is  isometrically isomorphic to $\ell_\infty$.
\end{co}
\pf Let $T$ be an injective, dense range, non surjective Tauberian operator on $L_1(0,1)$ (Theorem \ref{theorem2}).  Since $T$ is Tauberian, $T^{**}$ is also injective, so $T^*$ has dense  range but $T^*$ is not surjective because its range is not closed, and  $T^*$ is injective because $T$ has dense range.  The operator $T^*$ translates to an operator on $\ell_\infty$ that has the same properties because $L_\infty$ is isomorphic to $\ell_\infty$ by an old result due to Pe\l czy\'nski (see, e.g., \cite[Theorem 4.3.10]{ak}) (notice however that, unlike $T^*$,  the operator on $\ell_\infty$ cannot be weak$^*$ continuous).

For the ``consequently" statement, let $S$ be any norm one  injective, dense range, non surjective operator on $\ell_\infty$.  In the space $\ell_\infty \oplus_\infty \ell_\infty$, which is isometric to $\ell_\infty$, define
$X:= \ell_\infty \oplus \{0\}$ and $Y:= \{(x,Sx) : x\in \ell_\infty\}$. Obviously $X$ and $Y$ are isometric to $\ell_\infty$ and $X+Y= \ell_\infty \oplus S \ell_\infty$, which is a dense proper subspace of  $\ell_\infty \oplus_\infty \ell_\infty$.  Finally, $X\cap Y =\{0\}$ since $S$ is injective, so $X$ and $Y$ are quasi-complementary, non complementary subspaces of $\ell_\infty \oplus_\infty \ell_\infty$. \qed

Theorem \ref{theorem2} and the MathOverFlow question \cite{bah} suggest the following problem: Suppose $X$ is a separable Banach space (so that $X^*$ is isometric to a weak$^*$ closed subspace of $\ell_\infty$) and $X^*$ is non separable.  Is there a dense range operator on $X^*$ that is not surjective?  The answer is ``no":  Argyros, Arvanitakis, and Tolias \cite{aat} constructed a separable space $X$ so that $X^*$ is non separable, hereditarily indecomposable (HI), and every strictly singular operator on $X^*$ is weakly compact. Since $X^*$ is HI, every operator on $X^*$ is of the form $\lambda I +S$ with $S$ strictly singular.  If $\lambda \not= 0$, then  $\lambda I +S$  is Fredholm of index zero by Kato's classical perturbation theory. On the other hand, since every weakly compact  subset of the dual to a separable space is norm separable, every  strictly singular operator on $X^*$ has separable range.  (Thanks to Spiros Argyros for bringing this example to our attention.)

Any operator $T$ on $l^\infty$ that has dense range but is not surjective has the property that $0$ is an interior point of $\sigma(T)$. This follows from Thm 2.6 in \cite{A.N.}, where it is shown that $\partial\sigma(T)\subset\sigma_{p}(T^{*})$ for any operator $T$ acting on a $C(K)$ space which has the Grothendieck property.


 \begin{tabular}{l}
W.~B.~Johnson\\
Department of Mathematics\\
Texas A\&M University\\
College Station, TX  77843 U.S.A.\\
{\tt johnson@math.tamu.edu}
\\
\end{tabular}

\bigskip

  \begin{tabular}{l}
A.~B.~Nasseri\\
Fakult\"at f\"ur Mathematik\\
Technische Universit\"at Dortmund\\
 D-44221 Dortmund, Germany\\
{\tt amirbahman@hotmail.de}
\\
\end{tabular}

\bigskip

\begin{tabular}{l}
G.~Schechtman\\
Department of Mathematics\\
Weizmann Institute of Science\\
Rehovot, Israel\\
{\tt gideon@weizmann.ac.il}\\
\end{tabular}

\bigskip

\begin{tabular}{l}
T.~Tkocz\\
Mathematics Institute\\
University of Warwick\\
Coventry CV4 7AL, UK\\
{\tt ttkocz@gmail.com}
\end{tabular}

\end{document}